\documentclass[11pt]{amsart}
\usepackage{amsthm,amsmath,amssymb,amscd,epsfig}

{\end{list}}
{
   \newtheorem{theorem}{Theorem}[section]
   
   \newtheorem{lemma}[theorem]{Lemma}

   \newtheorem{conjecture}[theorem]{Conjecture}

}
{\theoremstyle{definition}

   \newtheorem{example}[theorem]{Example}
   \newtheorem{definition}[theorem]{Definition}
}
{\theoremstyle{remark}
   \newtheorem{remark}[theorem]{Remark}

}
\newcommand{\RR}{{\mathbb{R}}}

\newcommand{\CC}{{\mathbb{C}}}

\newcommand{\NN}{{\mathbb{N}}}

\newcommand{\cA}{{\mathcal A}}
\newcommand{\cB}{{\mathcal B}}
\newcommand{\cC}{{\mathcal C}}

\newcommand{\cE}{{\mathcal E}}
\newcommand{\cF}{{\mathcal F}}

\newcommand{\cL}{{\mathcal L}}

\newcommand{\cQ}{{\mathcal Q}}

\newcommand{\res}{\operatorname{res}}

\newcommand{\Star}{\operatorname{Star}}

\newcommand{\isom}{\simeq}

\newcommand{\Cd}{C^\bullet}
\newcommand{\Po}{Poincar\'e }

\begin{document}
\title{Lefschetz decomposition and the cd-index of fans.}

\author{Kalle Karu}
\thanks{The author was supported by NSERC}
\address{Mathematics Department\\ University of British Columbia \\
  1984 Mathematics Road\\
Vancouver, B.C. Canada V6T 1Z2}
\email{karu@math.ubc.ca}

\begin{abstract} 
The goal of this article is to give a Lefschetz type
  decomposition for the cd-index of a complete fan. 

To a complete simplicial fan one can associate a toric variety $X$,
  the even Betti 
  numbers $h_i$ of $X$ and the numbers $g_i = h_i-h_{i-1}$. If the fan
  is projective, then non-negativity of $g_i$ follows from the
  Lefschetz decomposition of the cohomology. 

In the case of a nonsiplicial complete fan one can analogously compute
the flag h-numbers $h_S$ and, by a change of variable formula, the cd-index. We
give an analogue of the Lefschetz operation for the cd-index. This
gives another proof of the non-negativity of the cd-index for complete fans. 
\end{abstract}

\maketitle

\section{Introduction}

Let $\Delta$ be a complete simplicial $n$-dimensional fan. Let $f_i$ be the number of $i$-dimensional cones in $\Delta$ and let $h_k$ be defined by the formula
\[ \sum_i f_{n-i}(t-1)^i = \sum_k h_{n-k} t^k.\] 
The numbers $h_k$ for $k=0,\ldots,n$ are the even Betti numbers $h_k =
\dim H^{2k}(X_\Delta,\CC)$ of a toric variety $X_\Delta$ if the fan
$\Delta$ is rational. If $\Delta$ is also projective, then there
exists a Lefschetz operation: 
\[ L:H^{2k}(X_\Delta,\CC) \to H^{2k+2}(X_\Delta,\CC), \qquad
L^k:H^{n-k}(X_\Delta,\CC) \stackrel{\isom}{\to}
H^{n+k}(X_\Delta,\CC),\] 
giving rise to the Lefschetz decomposition of the cohomology. The
existence of a Lefschetz operation implies that the numbers
$g_k=h_k-h_{k-1}$ are non-negative for $0\leq k\leq n/2$. 

For a complete but not necessarily simplicial fan one can construct
cohomology spaces $H^S(B\Delta)$ of dimension $h_S$ for $S\in \NN^n$,
and from the numbers $h_S$ compute the $cd$-index $\Psi_\Delta(c,d)$
of $\Delta$. Our goal is to find linear maps on $H^S(\Delta)$ that
guarantee non-negativity of the $cd$-index. Unlike the simplicial
case, it is not clear how such maps should be defined. We will give in
Definition~\ref{def-LO} a rather weak notion of a Lefschetz operation
which, nevertheless, is sufficient to imply non-negativity of the
$cd$-index. We also conjecture a stronger version in which the maps
are defined by conewise linear functions on the fan, just as in the
simplicial case. The rest of the introduction is spent on constructing
the $cd$-index and explaining the notion of a Lefschetz operation. 

Returning to the simplicial case, a simple way to construct the
cohomology $H^{2*}(X_\Delta,\CC)$ (which we will denote simply
$H^{*}(\Delta)$) is to consider the space $\cA(\Delta)$ of
complex-valued conewise polynomial functions on the fan $\Delta$. This
space is a free module under the action of the ring $A$ of global
polynomial functions, graded by degree. The graded vector space
$\cA(\Delta)/m\cA(\Delta)$, where $m\subset A$ is the maximal
homogeneous ideal, is the cohomology space
$H^{*}(\Delta)$ with \Po 
polynomial  
\[ P_\Delta(t) = \sum_k h_k t^k, \qquad h_k = \dim H^{k}(\Delta).\]
The fan $\Delta$ is projective iff there exists a strictly convex
conewise linear function $L\in \cA(\Delta)$. Multiplication with $L$
induces a Lefschetz operation in cohomology. 

In case when the fan $\Delta$ is complete, but not necessarily
simplicial, we proceed as follows. Let $B\Delta$ be a first
barycentric subdivision of $\Delta$. The space $\cA(B\Delta)$ is
graded by $\NN^n$, and with an adjustment of the module structure, the
quotient $H^*(B\Delta) := \cA(B\Delta)/ m\cA(B\Delta)$ inherits a
similar grading. Consider the corresponding \Po polynomial 
\[ P_{B\Delta}(t_1,\ldots,t_n) = \sum_{S\in \NN^n} h_S t_1^{S_1}\cdots
t_n^{S_n}, \qquad h_S = \dim H^S(B\Delta).\] 
\Po duality $h_S = h_{(1,\ldots,1)-S}$ implies that the sum can be
indexed by subsets $S\subset \{1,\ldots,n\}$.  

Let $\RR\langle c,d\rangle$ be the polynomial ring in non-commuting
variables $c$ and  $d$ of degree $1$ and $2$, respectively. There is
an embedding of vector spaces  
\[ \phi: \RR\langle c,d\rangle \hookrightarrow \RR[t_1,t_2,\ldots],\] 
defined as follows. $\phi$ maps constants to constants and if $f(c,d)c+g(c,d)d$ is a homogeneous $cd$-polynomial of degree $m>0$, define inductively 
\[ \phi(f(c,d)c+g(c,d)d) = \phi(f(c,d))(t_m+1)+\phi(g(c,d))(t_{m-1}+t_m).\]
For example, there are $3$ $cd$-monomials of degree $3$:
\begin{align*}
c^3 &= (t_1+1)(t_2+1)(t_3+1),\\
cd &= (t_1+1) (t_2+t_3),\\
dc &= (t_1+t_2) (t_3+1).
\end{align*}

It is shown in \cite{BK} that the \Po polynomial
$P_{B\Delta}(t_1,\ldots,t_n)$ of a complete fan $\Delta$ (more
generally, of a rank $n$ Eulerian poset) can be expressed as a
homogeneous $cd$-polynomial of degree $n$, called the $cd$-index
$\Psi_\Delta(c,d)$ of $\Delta$. The coefficients of the polynomial are
integers \cite{BK} and non-negative \cite{S,K}. 

\begin{figure}[ht] \label{fig-mon}
\centerline{\psfig{figure=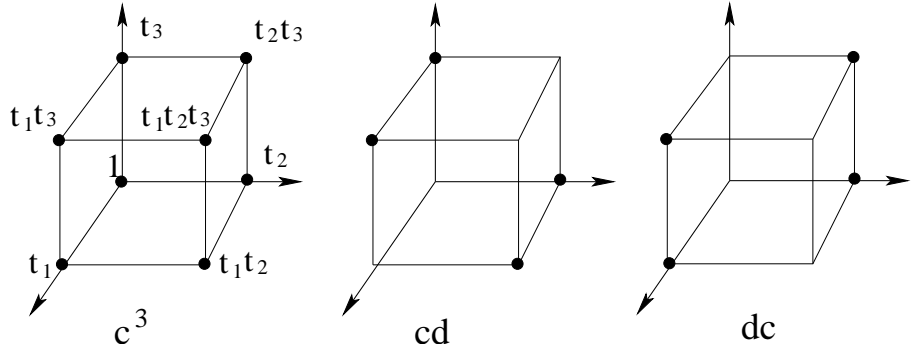,width=12cm}}
\caption{$cd$-monomials in terms of $t_i$-monomials.}
\end{figure}

One approach to proving non-negativity of the $cd$-index is to
decompose the cohomology $H^*(B\Delta)$ into summands corresponding to
different $cd$-monomials, such that the coefficients of
$\Psi_\Delta(c,d)$ are the dimensions of the corresponding
components. If we know the non-negativity of the $cd$-index, then the
existence of such a decomposition follows
trivially. Figure~\ref{fig-mon} shows the dimensions of the pieces
corresponding to different $cd$-monomials in the $3$-dimensional
case. The bold dots indicate the $t_i$-monomial being a summand of the
$cd$-monomial.  

In analogy with the singly-graded case we expect the
decomposition to be defined by linear maps. More precisely, we look
for endomorphisms $L_i: H^*(B\Delta) \to  H^*(B\Delta)$ of degree
$e_i=(0,\ldots,1,\ldots,0)$. If a $cd$-monomial $m$ can be written as
$m= \ldots(t_i+1)\ldots$, then $L_i$ should map in the corresponding
piece $H_m^*$ of the cohomology decomposition: 
\[ L_i: H_m^{(*,\ldots,*,0,*,\ldots,*)} 
\stackrel{\isom}{\longrightarrow} H_m^{(*,\ldots,*,1,*,\ldots,*)}.\] 
For example, $L_1$ should define an isomorphism from the back face of
the cube to the front face in Figure~\ref{fig-mon} for the monomials
$c^3$ and $cd$; the component corresponding to the monomial $dc$
should lie in the kernel of $L_1$.  

\subsection{The Main Construction} \label{sec-main-constr}

The definition of a Lefschetz operation is given inductively using a
construction that we call "the main construction". It essentially
describes the action of $L_1$ on the $A$-module $\cA(B\Delta)$ as
described in the paragraph above. 

Let $A_{l,m}$ be the polynomial ring $\CC[x_l,\ldots,x_m]$, graded by
$\NN^{m-l}$, with $x_i$ having degree $e_i$. Let the dualizing module
of $A_{l,m}$ be $\omega_{l,m}$, the principal ideal in $A_{l,m}$
generated by $x_l\cdots x_m$.

Let $M$ be a finitely generated free graded $A_{l,m}$-module. A \Po
pairing on $M$ is an $A_{l,m}$-bilinear symmetric map
\[ <\cdot,\cdot>: M\times M \to \omega_{l,m},\]
inducing a nondegenerate pairing on $\overline{M} =
M/(x_l,\ldots,x_m)M$. We always assume that $M$ is graded in
non-negative degrees. Then the existence of a \Po pairing implies that
$\overline{M}$ is graded by subsets of $\{l,\ldots,m\}$. 

Let $M$ be a free $A_{l,m}$-module with a \Po pairing and let $L:M\to
M$ be an endomorphism of degree $e_l$ which is self-adjoint with
respect to the pairing: 
\[ <L m_1, m_2> = <m_1, L m_2>.\]
We can write 
\[ M/(x_l)M = M^0 + M^1,\]
where $M^i$ consists of elements of degree $(i,*,\ldots,*)$. Then $L$
induces a map $M^0\to M^1$. 

Assume that the map  $L: M^0\to M^1$ is injective and the quotient is
annihilated by $x_{l+1}$: 
\begin{equation}\label{ass-main}
 0\to M^0 \stackrel{L}{\to} M^1 \to Q\to 0, \qquad x_{l+1} Q = 0.
\end{equation}
Then $Q$ is a free $A_{l+2,m}=\CC[x_{l+2},\ldots,x_m]$-module and we
get a long-exact Tor sequence: 
\begin{equation}\label{seq-A}
0\to Q[e_l-e_{l+1}]\to M^0/(x_{l+1})M^0 \stackrel{L}{\to}
M^1/(x_{l+1})M^1 \to Q\to 0.
\end{equation}
Let $C$ be the cokernel of the embedding $Q[e_l-e_{l+1}]\to M^0/(x_{l+1})M^0$:
\begin{equation}\label{seq-B}
0\to Q[e_l-e_{l+1}]\to M^0/(x_{l+1})M^0 \to C \to 0.
\end{equation}
Then $C$ is also a free $A_{l+2,m}$-module. We will show below that
$Q$ and $C$ both inherit a \Po pairing from $M$. The construction of
$Q$ and $C$ from $M$ and $L$ is what we call the main construction. 

Let $P_M(t_1,\ldots,t_n)$ (resp. $P_Q$, $P_C$) be the Hilbert
polynomial of $\overline{M}$ (resp. $\overline{Q}$,
$\overline{C}$). From the exact sequences (\ref{seq-A}) and (\ref{seq-B}),
we get 
\begin{equation} P_M = (1+t_1) P_C +(t_1+t_2)P_Q = c P_C +d
  P_Q.\label{seq-C}
\end{equation}
Thus, if $P_C$ and $P_Q$ are both $cd$-polynomials with non-negative
coefficients, then the same is true for $P_M$. We use this reason to
define Lefschetz operation inductively as follows. 

\begin{definition}\label{def-LO} Let $M$ be a finitely generated free
  $A_{l,m}$-module with a \Po pairing. We say that $M$ has a {\em Lefschetz
  operation} if there exists an endomorphism $L:M\to M$ of degree
  $e_l$, satisfying the assumptions of the main construction, such
  that the modules $C$ and $Q$ also have Lefschetz operations. More precisely: 
\begin{itemize}
\item $L$ is self-adjoint with respect to the pairing on $M$.
\item $L:M^0\to M^1$ is injective with cokernel annihilated by
  $x_{l+1}$.
\item Inductively, the $A_{l+1,m}$-module  $C' = C\otimes
  \CC[x_{l+1}]$ and the $A_{l+2,m}$-module $Q'=Q[e_l]$ have Lefschetz
  operations.
\end{itemize}
To start the induction, if $l>m$ and $M$ is a finite dimensional
vector space, then it trivially has a Lefschetz operation.
\end{definition}

Let us explain the replacement of $Q$ and $C$ by $Q'$ and $C'$,
respectively. Note that $Q$ lies in degrees $(1,0,*,\ldots,*)$, with
its \Po dual $Q[e_l-e_{l+1}]$ in degrees $(0,1,*,\ldots,*)$. Thus, to
get a degree preserving pairing on $Q$ to $\omega_{l+2,m}$, we have to
shift it by $e_l$. Going from $M$ to $Q$ and $C$ corresponds to the
$cd$-monomials $d$ and $c$, respectively. Since degree of $d$ is $2$
and the degree of $c$ is $1$, we expect $Q$ to be a $A_{l+2,m}$-module
with pairing into $\omega_{l+2,m}$ and $C$ to be a $A_{l+1,m}$-module
with pairing into $\omega_{l+1,m}$. Therefore we replace $C$ by the
$A_{l+1,m}$-module $C'$.

From the computation (\ref{seq-C}) above, it is clear that if  $M$ has a
Lefschetz operation, then the Hilbert function of $\overline{M}$  can
be written as a homogeneous $cd$-polynomial of degree $m-l$ with
non-negative integer coefficients. (To be precise, in the embedding
$\phi: \RR\langle c,d\rangle \hookrightarrow \RR[t_1,t_2,\ldots]$ we
need to relable the variables $t_i$ so that they sart with $t_l$.)

The main result of this article is:

\begin{theorem} \label{thm-main}
Let $\Delta$ be a complete fan of dimension $n$. Then the
$A_{1,n}$-module $\cA(B\Delta)$
has a Lefschetz operation. In particular, the $cd$-index of $\Delta$
has non-negative integer coefficients.  
\end{theorem}

Recall that $\cA(B\Delta)$ is a ring. If $L_i\in \cA(B\Delta)$
is an element of degree $e_i$, then multiplication with $L_i$ defines
an endomorphism of   $\cA(B\Delta)$ of degree $e_i$, self-adjoint with
respect to the natural \Po pairing. Thus, $L_1$ is a good candidate
for the Lefschetz operation on $\cA(B\Delta)$, and inductively, $L_i$
for $i>1$ could be used to define the endomorphisms of $Q$ and $C$.

\begin{conjecture} \label{conj-main} Let $L_i\in \cA(B\Delta)$ be a general element of
  degree $e_i$ for $i=1,\ldots,n$. Then $L_i$ define a Lefschetz
  operation on $\cA(B\Delta)$.
\end{conjecture}

We remark that a Lefschetz operation on $M$ does not define a
canonical decomposition of $\overline{M}$ into components
corresponding to the $cd$-monomials. To decompose $\overline{M}$, we
need to choose a splitting of the sequence~(\ref{seq-B}), so that 
\[ \overline{M} \isom \overline{C}\oplus \overline{C}[-e_l] \oplus
\overline{Q}\oplus \overline{Q}[e_l-e_{l+1}],\]
corresponding to the formula~(\ref{seq-C}). Inductive decomposition of
$\overline{C}$ and  $\overline{Q}$ then give a complete decomposition
of $\overline{M}$.

To prove Theorem~\ref{thm-main}, we express $\cA(B\Delta)$ as the space
of global sections of a sheaf $\cL$ on $\Delta$. The main construction
can be sheafified, i.e., performed on the stalks of the sheaf $\cL$
simultaneously. We show that a Lefschetz operation on the space of
global sections comes from a sheaf  homomorphism.

We also consider Conjecture~\ref{conj-main} in the context of sheaves
and reduce it to a Kleiman-Bertini type problem of torus actions on a vector
space. Let an algebraic torus $T$ act on a finite dimensional vector
space $V$ with possibly infinitely many orbits. When does the general
translate of a subspace $K\subset V$ intersect another subspace
transversely? Conjecture~\ref{conj-last} states sufficient
conditions for this, implying  Conjecture~\ref{conj-main}. 

Theorem~\ref{thm-main} gives another proof of non-negativity of the
$cd$-index for a complete fan. In \cite{K} 
non-negativity was proved more generally for Gorenstein* posets. The
current proof does not extend to that more general situation. The two
proofs are based on the same idea. However,
the proof we give here is simpler because we work with modules only,
avoiding derived categories.

\section{Sheaves on Fans}

All our vector spaces are over the field of complex numbers $\CC$. 
Let $A_{l,m}=\CC[x_l,x_{l+1},\ldots,l_m]$, graded so that $x_i$ has
degree $e_i$.  For a graded $A_{l,m}$-module $M$ we denote the shift
in grading by $M[\cdot]$. We also write $\overline{M} = M/(x_l,\ldots,x_m)M$. 

For a graded set, the superscript refers to degree. If $\Delta$ is a
fan, then $\Delta^{\geq m}$ consists of all cones of dimension at
least $m$. Similarly, $\Delta^{[l,m]}$ is the subset of cones of
dimension $d\in[l,m]$.

\subsection{Fan spaces} 
Let us recall the notion of sheaves on fans. The main reference for
the general theory is \cite{BBFK,BL} and for the specific sheaves used
here \cite{K}. 

We fix a complete $n$-dimensional fan $\Delta$ (see \cite{F} for
terminology). Consider $\Delta$ as a finite partially ordered set of
cones, graded in degrees $0,\ldots,n$. It is sometimes convenient to
add a maximal element $\hat{1}$ of degree $n+1$ to $\Delta$.  

The fan $\Delta$ is given the topology in which open sets are the
(closed) subfans of $\Delta$. Then a sheaf $F$ of vector spaces on
$\Delta$ consists of the data: 
\begin{itemize}
\item A vector space $F_\sigma$ for each $\sigma\in\Delta$.
\item Linear maps $res^\sigma_\tau: F_\sigma \to F_\tau$ for $\sigma>\tau$, satisfying the compatibility condition $res^\tau_\rho \circ res^\sigma_\tau = res^\sigma_\rho$ for $\sigma>\tau>\rho$. 
\end{itemize}

On sheaves we can perform the usual sheaf operations. For example, a
global section $f\in\Gamma(F,\Delta)$ consists of the data
$f_\sigma\in F_\sigma$ for each $\sigma\in\Delta$, such that
$res^\sigma_\tau f_\sigma = f_\tau$. Equivalently, we only need to
give $f_\sigma\in F_\sigma$ for maximal cones $\sigma$, such that their
restrictions to smaller dimensional cones agree.

Define a sheaf of rings $\cA$ on $\Delta$ as follows:
\begin{itemize}
\item $\cA_\sigma = A_{1,d}= \CC[x_1,\ldots,x_d]$ if $\dim\sigma = d$.
\item $res^\sigma_\tau: \CC[x_1,\ldots,x_d] \to \CC[x_1,\ldots,x_l]$ is
  the standard projection $x_i\mapsto x_i$ for $1\leq i \leq l$ and $x_i\mapsto
  0$ for $i>l$. 
\end{itemize}

Given the sheaf of rings $\cA$ on $\Delta$, we consider sheaves of
$\cA$-modules $\cF$. This means that the stalks $\cF_\sigma$ are
$\cA_\sigma$-modules and the restriction maps are module
homomorphisms. Note that the sheaf $\cA$ is graded by $\NN^n$. We
assume that all sheaves of $\cA$-modules are similarly graded. 

There exists an indecomposable sheaf $\cL$ of $\cA$-modules satisfying
the following conditions: 
\begin{itemize}
\item Locally free: $\cL_\sigma$ is a graded free $\cA_\sigma$-module.
\item Minimally flabby: $\dim \cL_0 = 1$ and for $\sigma>0$, the
  restriction maps induce an isomorphism 
\[ \overline{\cL}_\sigma \to \overline{\Gamma(\cL,\partial\sigma)},\]
where $\partial\sigma$ is the boundary fan of $\sigma$.
\end{itemize}

These two conditions define $\cL$ up to an isomorphism. In fact,
$\Gamma(\cL,\partial\sigma)$ is a free
$A_{1,d-1}$-module if $\dim\sigma = d$, and we can
inductively define 
\[ \cL_\sigma = \Gamma(\cL,\partial\sigma) \otimes_{A_{1,d-1}} A_{1,d}.\]

\subsection{Barycentric Subdivisions}

Let $B\Delta$ be a barycentric subdivision of $\Delta$. As a poset it consists of chains $x=(0<\sigma_1<\ldots<\sigma_m)$ in $\Delta$. 
Define a sheaf of rings $\cB$ on $B\Delta$ as follows:
\begin{itemize}
\item $\cB_x = \CC[x_i]_{i\in S}$, where $x=(0<\sigma_1<\ldots<\sigma_m)$, $S = \{ \dim \sigma_1,\ldots,\dim\sigma_m\}$.
\item $res^x_y$ is the standard projection.
\end{itemize}

One can construct as above a sheaf $\cL$ with respect to $\cB$, but this sheaf is isomorphic to $\cB$. 

\begin{lemma}[\cite{K}] We have
\[ \pi_* \cB \isom \cL,\]
where $\pi: B\Delta\to \Delta$ is the subdivision map sending
$x=(0<\sigma_1<\ldots<\sigma_m)$ to $\sigma_m$. \qed
\end{lemma}

It is often more convenient to work with the sheaf $\cB$ because it is
a sheaf of rings. The space of global sections $\Gamma(\cB,B\Delta)$
(which is isomorphic to $\Gamma(\cL,\Delta)$ by the previous lemma) is
what we called $\cA(B\Delta)$ in the introduction. Since $\cB$ and
$\cL$ are sheaves of $A_{1,n}$-modules, so are the spaces of global
sections. 

\subsection{The Cellular Complex}

Let us fix an orientation for each cone $\sigma\in\Delta$ and for $\sigma>\tau$, $\dim\sigma=\dim\tau+1$, let 
\[ or^\sigma_\tau = \pm 1\]
depending on whether the orientations of $\sigma$ and $\tau$ agree or not.

The cellular complex of a sheaf $F$ on $\Delta$ is 
\[ \Cd_n(F,\Delta) = 0\to C^0 \to C^1 \to \ldots \to C^n\to 0,\]
where
\[ C^i = \bigoplus_{\dim\sigma = n-i} F_\sigma,\]
and the differentials are defined as sums of $or^\sigma_\tau res^\sigma_\tau:
F_\sigma\to F_\tau$. 

For a complete fan $\Delta$, the cellular complex $\Cd_n(F,\Delta)$
computes the cohomology of $F$. Applying this to the flabby sheaf
$\cL$, we get 
\[ H^i(\Cd_n(\cL,\Delta)) = \begin{cases} \Gamma(\cL,\Delta) & \text{if $i=0$} \\
					0	& \text{otherwise.}
                           \end{cases}
\]
Moreover, $\Gamma(\cL,\Delta)$ is a graded free $A_{1,n}$-module.

If $\sigma\in\Delta$, $\dim\sigma = d$, then $\partial\sigma$ is
combinatorially equivalent to a complete fan of dimension $d-1$, hence
we may use $\Cd_{d-1}(\cL,\partial\sigma)$ to compute
$\Gamma(\cL,\partial\sigma)$. This gives an exact sequence 
\[ 0\to L_\sigma/x_d L_\sigma \to \bigoplus_{\tau<\sigma,
  \dim\tau=d-1} \cL_\tau \to \bigoplus_{\rho<\sigma, \dim\rho=d-2}
\cL_\rho \to \ldots\to \cL_0 \to 0.\] 

\subsection{\Po Pairing}

Define the dualizing module $\omega_{1,n} = (x_1\cdots
x_n)A_{1,n}$. I.e., $\omega_{1,n}$ is the principal ideal generated by
$x_1\cdots x_n$.  There exists an $A_{1,n}$-bilinear non-degenerate pairing  
\[ \Gamma(\cL,\Delta) \times \Gamma(\cL,\Delta) \to \omega_{1,n}.\]

The pairing is best constructed using the isomorphism $\Gamma(\cL,\Delta)
\isom \Gamma(\cB,B\Delta)$. On $\Gamma(\cB,B\Delta)$ the pairing is
defined by multiplication ($\cB$ is a sheaf of rings), followed by an
evaluation map into $\omega_{1,n}$.  

One can give a simple description of the evaluation map as in
\cite{B}, depending on the orientations $or^\sigma_\tau$. For
$x=(0<\sigma_1<\ldots<\sigma_n)$ a maximal element of $B\Delta$ of
dimension $n$, define  
\[ \varepsilon_x = or^{\hat{1}}_{\sigma_n}
or^{\sigma_n}_{\sigma_{n-1}}\cdots or^{\sigma_1}_0 = \pm 1.\] 
Now if $f\in \Gamma(\cB,B\Delta)$, then it can be shown that
\[ \sum_{\dim{x}=n} \varepsilon_x f_x\]
is an element of $A_n$ that is divisible by $x_1 x_2\cdots x_n$, hence
lies in $\omega_n$. This defines the
$A_n$-linear evaluation map $\Gamma(\cB,B\Delta)\to \omega_n$ and the
\Po pairing on $\Gamma(\cB,B\Delta)$. 

If $\sigma\in\Delta$ is a $d$-dimensional cone, then $\partial\sigma$
is combinatorially equivalent to a complete fan of dimension
$d-1$. By the same construction as above we get a pairing on $\Gamma(\cL,
\partial\sigma) \isom \cL_\sigma/ x_d\cL_\sigma$. 

In summary, for each cone $\sigma\in\Delta$, $\dim\sigma=d$, we have a
non-degenerate symmetric bilinear pairing 
\[ <\cdot,\cdot>_\sigma: \cL_\sigma/ x_d\cL_\sigma \times \cL_\sigma/
x_d\cL_\sigma \to \omega_{1,d-1}.\]
These pairings are related as follows. For $f,g\in \cL_\sigma/
x_d\cL_\sigma$, 
\[ <f,g>_\sigma = \sum_{\dim\tau=d-1} or^\sigma_\tau
<f_\tau,g_\tau>_\tau,\]
where $f_\tau$ and $g_\tau$ are the restrictions of $f$ and $g$ to
$\tau$ and the pairing on the right hand side is the
$A_{1,d-1}$-bilinear extension of the $A_{1,d-2}$-bilinear pairing
$<\cdot,\cdot>_\tau$.

\section{The Main Construction on Sheaves}

Let us return to the situation of Section~\ref{sec-main-constr} and prove the claims made there. 

We have a finitely generated free $A_{l,m}$-module $M$ with \Po pairing
\[ <\cdot,\cdot>_M: M\times M \to \omega_{l,m}.\]
Write
\[ M/x_l M = M^0\oplus M^1,\]
where $M^i$ consists of elements of degree $(i,*,\ldots,*)$. Assume
that $L:M^0\to M^1$ is a $A_{l+1,m}$-module homomorphism of degree
$e_l$, self-adjoint with respect to the pairing, and such that $L$ is
injective with quotient $Q$ annihilated by $x_{l+1}$: 
 \[ 0\to M^0 \stackrel{L}{\to} M^1 \to Q\to 0, \qquad x_{l+1} Q = 0. \label{ass-main1}\]

\begin{lemma} \label{lem-free}
Q is a free $A_{l+2,m}$-module.
\end{lemma}

{\bf Proof.} Since $M^0$ and $M^1$ are free $A_{l+1,m}$-modules, we
get from the exact sequence above that  
\[ Tor^{A_{l+1,m}}_i(Q,\CC) = 0, \qquad i\geq 2.\]
Because $Q$ is a $A_{l+2,m}$-module, annihilated by $X_{l+1}$, this
implies that 
\[ Tor^{A_{l+2,m}}_1(Q,\CC) = 0,\]
hence $Q$ is free. \qed

Now assuming that $Q$ is free, we get an exact sequence
\[0\to Q[e_l-e_{l+1}]\to M^0/(x_{l+1})M^0 \stackrel{L}{\to}
M^1/(x_{l+1})M^1 \to Q\to 0, \label{seq-A1}\] 
where all terms are free $A_{l+2,m}$-modules. Define $C$ by the exact
sequence 
\[0\to Q[e_l-e_{l+1}]\to M^0/(x_{l+1})M^0 \to C\to 0.\]
Then $Tor^{A_{l+2,m}}_1(C,\CC) = 0$ and $C$ is also a free $A_{l+2,m}$-module.

Let us construct bilinear pairings on $C$ and $Q$. On $C$ the pairing is 
\[ <x,y>_C = <x,Ly>_M.\]
This is well-defined and gives an $A_{l+2,m}$-linear map of degree $e_l$
\[ C\otimes_{A_{l+2,m}} C \to \omega_{l,m}\otimes_{A_{l,m}} A_{l+2,m}.\]
Dividing by $x_l$ we get a degree $0$ map into
$\omega_{l+1,m}\otimes_{A_{l+1,m}} A_{l+2,m}$. Finally, replacing $C$
by $C'= C \otimes_{A_{l+2,m}} A_{l+1,n}$ and extending the pairing
linearly, we have a $A_{l+1,m}$-bilinear map  
\[ <\cdot,\cdot>_{C'}: C'\times C' \to \omega_{l+1,m}.\]

To define the pairing on $Q$, let $\alpha$ be the composition
\[\alpha: Q\stackrel{\isom}{\to} Q[e_l-e_{l+1}] \hookrightarrow
M^0/x_{l+1} M^0.\] 
On the elements $[q]\in Q$ this map is given by
\[ \alpha([q]) = L^{-1}(x_{l+1} q).\]
Now define the pairing
\[ <x,y>_Q = <\alpha(x),y>_M.\]
One can  check that this pairing is well-defined. Since $Q$ lies
in degrees $(1,0,\ldots,*)$, we replace it with $Q' =
Q[e_l]$. Then, taking into account that $\alpha$ has degree $e_{l+1}-e_l$,
we get a degree $0$ $A_{l+2,m}$-bilinear map 
\[ <\cdot,\cdot>_{Q'}: Q'\times Q' \to \omega_{l+2,m}.\]
 
It is easy to see that the bilinear maps on $C'$ and $Q'$ are symmetric. 

\begin{lemma}
The pairings $<\cdot,\cdot>_{Q'}$ and $<\cdot,\cdot>_{C'}$ are nondegenerate.
\end{lemma}

{\bf Proof.} One  checks the non-degeneracy of the pairing on
$C$ using the definition and self-adjointness of $L$. Then it follows
that the pairing between $Q$ and $Q[e_l-e_{l+1}]$ is
non-degenerate. \qed

We next want to sheafify the main construction. Recall that $\cL$ is a
sheaf on $\Delta$ with stalks $\cL_\sigma$ free $\cA_\sigma$-modules with \Po
pairings. To perform the main construction simultaneously on all 
stalks of $\cL$, the first step is to split 
\[ \cL/x_1\cL = \cL^0\oplus \cL^1,\]
and then find a map of sheaves of degree $e_1$
\[ L: \cL^0 \to \cL^1.\]
If one looks at the stalks, it becomes clear that $\cL^i$ should be
considered as sheaves on $\Delta^{\geq 2}$ (i.e., on the poset of
cones of dimension at least $2$), and the cokernel $Q$ of the map $L$
should be a sheaf on $\Delta^{\geq 3}$. Therefore we will
consider sheaves on $\Delta^{\geq m}$ for $m\geq 1$.

\subsection{Sheaves on $\Delta^{\geq m}$}

We let $\Delta^{\geq m}$ have the the topology induced from
$\Delta$. To give a sheaf on 
$\Delta^{\geq m}$ is equivalent to giving a sheaf on $\Delta$ with all
stalks zero on cones of dimension less than $m$.

Define the structure sheaf $\cA$ on $\Delta^{\geq m}$ as follows. For $\sigma\in\Delta$, $\dim\sigma=d\geq m$, let
\[ \cA_\sigma = A_{m,d} = \CC[x_m,\ldots,x_d],\]
with restriction maps $res^\sigma_\tau$ the standard projections.

\begin{definition} Let $\cF$ be a locally free sheaf of $\cA$-modules
  on $\Delta^{\geq m}$. We say that $\cF$ is {\em minimally flabby} if
  all the restriction maps $res^\alpha_\beta$ are surjective and for
  every $\sigma\in\Delta$, $\dim\sigma=d\geq m$, we have an exact
  sequence, the "augmented cellular complex" 
\begin{equation} \label{seq-D} 0\to \cF_\sigma/x_d \cF_\sigma \to
  \bigoplus_{\tau<\sigma, \dim\tau=d-1} \cF_\tau \to \ldots \to
  \bigoplus_{\rho<\sigma, \dim\rho=m} \cF_\rho \to G_\sigma \to 0,
\end{equation} 
where
\begin{itemize}
\item The augmentation $G_\sigma$ is a vector space (i.e., an
  $A_{1,n}$-module annihilated by $x_1,\ldots,x_n$).
\item The differentials are defined by $or^\alpha_\beta res^\alpha_\beta$ as in the usual cellular complex.
\end{itemize}
\end{definition}

\begin{remark} \begin{enumerate}
\item It should be noted that a minimally flabby sheaf is not flabby
  in the topology of $\Delta^{\geq m}$.  
\item We do not need the surjectivity of the restriction maps
  $res^\alpha_\beta$ for the proof of Theorem~\ref{thm-main}. These
  conditions are only necessary to state Conjectures \ref{conj-main}
  and \ref{conj-last}. However, surjectivity of the restriction maps
  follows easily for all sheaves we consider.
\end{enumerate}
\end{remark}

\begin{example}
\begin{enumerate}
\item Let $\cL$ be the indecomposable sheaf on $\Delta$. Then
  $\cL|_{\Delta^{\geq 1}}$ is a minimally flabby sheaf on
  $\Delta^{\geq 1}$. In this case we have $G_\sigma = \cL_0 = \CC$ for
  all $\sigma$. 
\item In general, the vector spaces $G_\sigma$ depend on the cone
  $\sigma$. Let $\pi_1, \pi_2\in\Delta$ be two cones of dimension $m-1$, and
  let $\cL^{\pi_i}$ be the indecomposable sheaf constructed on the
  poset $\Star \pi_i$. Then $\cF = \cL^{\pi_1} \oplus \cL^{\pi_2}
  |_{\Delta^{\geq m}}$ is a minimally flabby sheaf and we have  
\[ G_\sigma = \begin{cases} 
\CC \oplus \CC & \text{if $\pi_1,\pi_2 <\sigma$} \\
\CC & \text{if $\pi_1<\sigma$ or $\pi_2 <\sigma$, but not both} \\
0 & \text{otherwise}.
              \end{cases}
\]
\end{enumerate}
\end{example}

Note that a minimally flabby sheaf on $\Delta^{\geq m}$ is determined
by its restriction to $\Delta^{[m,m+1]}$. Indeed, the exact sequence~(\ref{seq-D}) can be used to recover
$\cF_\sigma$ for $\dim\sigma > m+1$. Similarly, given two minimally
flabby shaves $\cF$ and $\cE$, a morphism defined between the
restrictions of these sheaves to $\Delta^{[m,m+1]}$ can be lifted to a
morphism on $\Delta^{\geq m}$.

\begin{lemma}\label{lem-quot} Let $\cE$ and $\cF$ be minimally flabby
  sheaves on $\Delta^{\geq m}$, and $L: \cE\to \cF$ a homomorphism of
  $\cA$-modules.  
\begin{enumerate}
\item If $L$ is injective on cones $\sigma\in\Delta$, $\dim\sigma =
  m$, then $L$ is injective on all cones. 
\item If $L$ is an isomorphism on cones $\sigma\in\Delta$, $\dim\sigma
  = m$, then the cokernel $\cQ$ of $L$: 
\[ 0\to \cE\to \cF \to \cQ\to 0\]
is a minimally flabby sheaf on $\Delta^{\geq m+1}$.
\end{enumerate}
\end{lemma}

{\bf Proof.} The first statement follows by induction on $\dim\sigma$
from the exact sequence~(\ref{seq-D}). 

To prove the second statement, first note that the surjectivity of
the restriction maps $res^\alpha_\beta$ for $\cQ$ is clear.
The morphism $L$ defines a map between the augmented cellular
complexes of $\cE$ and $\cF$ which is injective except possibly in the
$G_\sigma$ terms. The quotient gives the cellular complex for $Q$. By
induction on $\dim\sigma$ it follows that $Q_\sigma$ is annihilated by
$x_m$, hence is a free $A_{m+1,d}$-module by
Lemma~\ref{lem-free}. We get the augmentation for $\cQ$ by removing
tha augmentations of $\cE$ and $\cF$ and considering the long-exact
cohomology sequence of the short-exact sequence of complexes.
\qed

\begin{definition} Let $\cF$ be a minimally flabby sheaf on
  $\Delta^{\geq m}$. We say that $\cF$ is a \Po sheaf if for every
  $\sigma\in\Delta$, $\dim\sigma=d\geq m$, we have an
  $A_{m,d-1}$-bilinear non-degenerate symmetric pairing 
\[ <\cdot,\cdot>_\sigma: F_\sigma/x_d F_\sigma \times  F_\sigma/x_d
  F_\sigma \to \omega_{m,d-1},\] 
satisfying the compatibility condition:
\begin{equation} \label{eq-poinc}
 <f,g>_\sigma = \sum_{\tau < \sigma} or^\sigma_\tau
<\res^\sigma_\tau f, \res^\sigma_\tau  g>_\tau, \quad f,g\in
F_\sigma/x_d F_\sigma.
\end{equation}
Here on the right hand side $<\cdot,\cdot>_\tau$ denotes the
$A_{m,d-1}$-bilinear extension of the $A_{m,d-2}$-bilinear pairing
$<\cdot,\cdot>_\tau$.
\end{definition}

\begin{example} 
The sheaf $\cL|_{\Delta^{\geq 1}}$ is a \Po sheaf on $\Delta^{\geq 1}$.
\end{example}

Let $\cF$ be a \Po sheaf on $\Delta^{\geq m}$. Then
$\overline{\cF}_\sigma$ for $\dim\sigma = d\geq m$ is a vector space
graded by subsets of $\{m,\ldots,d-1\}$.  
Write $\cF/x_m\cF$ for the sheaf with stalks
\[ (\cF/x_m\cF)_\sigma = \cF_\sigma/x_m\cF_\sigma.\]
This is a locally free sheaf on $\Delta^{\geq m+1}$, and we can split it as
\[ \cF/x_m\cF = \cF^0\oplus \cF^1,\]
where $\cF^i_\sigma$ consists of elements of degree $(i,*,\ldots,*)$.

\begin{lemma}
Let $\cF$ be a \Po sheaf on $\Delta^{\geq m}$. Then $\cF^0$ and $\cF^1$ are minimally flabby sheaves on $\Delta^{\geq m+1}$.
\end{lemma}

{\bf Proof.} Let us cut the sequence~(\ref{seq-D}) into two exact sequences
\begin{gather*} 0\to \cF_\sigma/x_d \cF_\sigma \to
  \bigoplus_{\tau<\sigma, \dim\tau=d-1} \cF_\tau \to \ldots \to S \to 0, \\
0\to S\to \bigoplus_{\rho<\sigma, \dim\rho=m} \cF_\rho \to G_\sigma \to 0.
\end{gather*}
From the second sequence we get that $S$ is a free $\CC[x_m]$-module,
hence the first sequence remains exact after taking quotient by the
ideal $(x_m)$ and splitting into two according to degree. The two
sequences are the augmented cellular complexes for   $\cF^0$ and $\cF^1$.
\qed

Now we are ready to define the sheafified version of the main
construction. Let $\cF$ be a \Po sheaf on $\Delta^{\geq m}$ and
$L:\cF\to \cF$ an endomorphism of $\cA$-modules of degree $e_m$, such
that $L_\sigma:\cF_\sigma\to \cF_\sigma$ is self-adjoint with respect
to the pairing for each $\sigma$. (More precisely, $L_\sigma:
\cF_\sigma \to \cF_\sigma$ has to be self-adjoint with respect to the
$A_{m,d}$-linear extension of the pairing $<\cdot,\cdot>_\sigma$.)
Assume that the induced morphism $L: \cF^0\to \cF^1$ is injective on cones
$\sigma\in\Delta$, $\dim\sigma=m+1$; then it is an isomorphism on
these cones by \Po duality. Lemma~\ref{lem-quot} gives an exact
sequence 
\[ 0\to \cF^0\to \cF^1 \to \cQ\to 0,\]
where $\cQ$ is a minimally flabby sheaf on $\Delta^{\geq m+2}$ In
order to have $\cQ$ in correct degrees, we have to replace it with
$\cQ' = \cQ[e_{m}]$.  

We also construct the sheaf $\cC$ as follows. First, we have an exact
sequence of minimally flabby sheaves on $\Delta^{\geq m+2}$: 
\[ 0\to \cQ[e_{m}-e_{m+1}]\to  \cF^0/x_{m+1}\cF^0 \to
\cF^1/x_{m+1}\cF^1 \to \cQ\to 0.\] 
Define $\cC$ by the exact sequence
\[ 0\to \cQ[e_{m}-e_{m+1}]\to  \cF^0/x_{m+1}\cF^0 \to \cC \to 0.\]
Then one easily sees that $\cC$ is also minimally flabby on
$\Delta^{\geq m+2}$ (to get the augmented cellular complex for $\cC$,
it is more convenient to consider the short exact sequence 
\[ 0\to \cC \stackrel{L}{\to} \cF^1/x_{m+1}\cF^1 \to \cQ\to 0).\]
We should again replace $\cC$ with an almost
flabby sheaf $\cC'$ on  $\Delta^{\geq m+1}$, such that $\cC =
\cC'/x_{m+1}\cC'$. We will not do this because inductively, the next
step to construct a Lefschetz operation is to go from $\cC'$ to $\cC$
and split it according to degree. The fact that we don't have $\cC'$
that induces $\cC$ will cause us some trouble later when we look for
an endomorphism of $\cC$. 

Summarizing, we have defined the sheafified version of the main
construction. Starting with a \Po sheaf $\cF$ on $\Delta^{\geq m}$ and
a morphism $L$, we constructed minimally flabby sheaves $\cQ$ and
$\cC$ on $\Delta^{\geq m+2}$. The construction on stalks agrees with
the main construction on modules. The stalks of the shaves $\cQ$ and
$\cC$ inherit \Po pairings from the pairing on $\cF$, which is clearly
compatible with the restriction morphisms. Hence the two new sheaves
are also \Po sheaves. 

It remains to see when can we find an appropriate endomorphism $L$ of $\cF$.

\begin{lemma} Let  $\cF$ be a \Po sheaf on $\Delta^{\geq m}$ and $L:
  \cF \to \cF$ a homomorphism of degree $e_m$. Then $L$ is
  self-adjoint with respect to the pairings on $\sigma\in
  \Delta^{\geq m}$ if and only if it is self-adjoint on cones $\rho$
  of dimension $m$.
\end{lemma}

{\bf Proof.} This follows by induction on the dimension of a cone from
the formula~(\ref{eq-poinc}). \qed

\begin{lemma}\label{lem-exist}
Let $\cF$ be a \Po sheaf on $\Delta^{\geq m}$. Then there exists a
homomorphism $L: \cF\to \cF$ of degree $e_m$ that is
self-adjoint with respect to the pairings on the stalks $\cF_\sigma$
and such that the induced homomorphism $L: \cF^0\to \cF^1$ is
injective.
\end{lemma}

{\bf Proof.} For $\dim\rho = m$, let $L_\rho: \cF_\rho \to \cF_\rho$
be a self-adjoint homomorphism of degree $e_m$. (Note that $\cF_\rho \isom
\CC[x_m]^{\oplus a_\rho}$ for some $a_\rho \geq 0$.) We claim that a
suitable collection of $L_\rho$ induces the required $L$. For this we
need to check that $L_\rho$ can be extended to cones $\tau$ of
dimension $m+1$ (hence can be extended to all cones), and that on such
$\tau$ it defines an injection $\cF_\tau^0\to \cF_\tau^1$.

Let $\dim\tau=m+1$ and consider the augmented cellular complex of
$\tau$: 
\[ 0\to \cF_\tau/x_{m+1}\cF_\tau  \to  \bigoplus_{\rho<\tau} \cF_\rho
\to G_\tau\to 0.\] 
Here $G_\tau\isom \CC^a$ for some $a\geq 0$, $\bigoplus_{\rho<\tau}\cF_\rho \isom \CC[x_m]^{\oplus 2a}$ and $\cF_\tau/x_{m+1}\cF_\tau \isom \CC[x_m]^{\oplus a} \oplus \CC[x_m][-e_m]^{\oplus a}$.

The maps $L_\rho$ are compatible with the zero map $G_\tau\to G_\tau$
of the augmentation. It follows that $L_\rho$ induce a map $L_\tau:
\cF_\tau \to \cF_\tau$, compatible with restriction maps, hence there
is an extension to a morphism $L:\cF\to\cF$.  

Let $V = \bigoplus_{\rho<\tau} \overline{\cF}_\rho \isom \CC^{2a}$ and let
$K\isom \CC^a$ be the kernel of $V \to G_\tau$. Then $K =
\cF^0_\tau$. The map $L_\rho$ comes from a linear map
$\overline{L}_\rho = \frac{1}{x} L_\rho:
\overline{\cF}_\rho \to \overline{\cF}_\rho$. The maps
$\overline{L}_\rho$ together define a linear map $L_V:V\to V$.  Now
the condition that $L$ is injective is equivalent to $L_\tau:
\cF_\tau^0\to \cF_\tau^1$ 
being injective, which is equivalent to the condition that the
intersection of $K$ and $L_V(K)$ is zero. 

Let us also bring the \Po pairing into the picture. We have a
non-degenerate symmetric pairing on each $\overline{\cF}_\rho$, combined to a pairing on $V$. The pairing on $\cF_\tau$ induces a
non-degenerate pairing between $\cF_\tau^0$ and $\cF_\tau^1$, which
restricts to the zero pairing on $\cF_\tau^0$, hence the compatibility
condition implies that the pairing on $V$ restricted to $K$ is
zero. In other words, $K=K^\perp$. The proof that a suitable set of
$L_\rho$ gives a required $L$ is given in the lemma below. 

Finally, let us consider the case when $L$ is defined by a multiplication
with an element in $L\in \Gamma(\cA,\Delta)$ of degree $e_m$. In this
case the linear maps $\overline{L}_\rho$ are given my multiplication with a
constant $c_\rho$ (where $L|_\rho = c_\rho x_m$). Note also that since
the restriction maps $res^\tau_\rho$ are surjective, the projection
$V\to \overline{\cF}_\rho$ maps $K$ onto $\overline{\cF}_\rho$.
Thus, if the conjecture below is true then $L$ defines an injective
morphism. \qed  

\begin{lemma}\label{lem-last} Let $V = \oplus V_i$ be a finite
    dimensional vector 
  space. Suppose that each $V_i$ has a non-degenerate symmetric
  bilinear pairing, giving a pairing on $V$. Let $K\subset V$ be a subspace
  such that $K \subset K^\perp$. Then there exist self-adjoint linear
  maps $L_i: V_i\to V_i$, combined 
 to $L:V\to V$, satisfying $K^\perp \cap L(K) = 0$. 
\end{lemma}

{\bf Proof.}
Let $v_1,\ldots,v_{2a}$ be an orthogonal basis of $V$ consisting of
elements from $V_i$ and let $y_1,\ldots,y_{2a}$ be the dual basis
giving coordinates on $V$. Let $T$ be the algebraic
torus of dimension $\dim V$ acting on $V$ by: 
\[ (t_1,\ldots,t_{2a})\cdot (y_1,\ldots,y_{2a}) = (t_1
y_1,\ldots,t_{2a}y_{2a}).\] 
An element $t\in T$ defines a linear map $V\to V$ of the required
type. We claim that for a general $t$ we have $K^\perp\cap t(K) = 0$. 
  
Now $V$ has finitely many $T$-orbits. By Kleiman-Bertini theorem, for
a general $t$, the restrictions of $K^\perp$ and $K$ to any orbit $O$
intersect transversely. Thus, it suffices to show that the expected
dimension of this intersection is zero.

Let $W\subset V$ be a subspace spanned by a subset of the $v_j$. Then
the pairing on $V$ restricts to a non-degenerate pairing on $W$. Since
$K \subset K^\perp$, it follows that 
\[ \dim(K^\perp\cap W) + \dim (K\cap W) \leq \dim(W). \qed \]

\begin{conjecture}\label{conj-last} Let the notation be as in the
  previous lemma. Additionally assume that the projections $V\to V_i$
  map $K$ onto $V_i$ for each $i$.Then the statement of the lemma
  remains true if we let $L_i$ be multiplication by some constant $c_i$. 
\end{conjecture}

\begin{remark}
Starting with a \Po sheaf $\cF$ on $\Delta^{\geq m}$, we apply the
previous lemmas to perform the main construction on $\cF$ and produce
new sheaves $\cQ$ and $\cC$. Then inductively we apply the same
construction on $\cC$ and $\cQ$. As explained above, we should
consider $\cC$ as coming from a sheaf $\cC'$ on $\Delta^{\geq m+1}$,
so that the main construction should be applied to $\cC'$ rather than
$\cC$. Let us show that we don't need $\cC'$ for the existence of the
required $L:\cC\to \cC$.

Recall that $\cC$ was defined by the exact sequence of minimally flabby sheaves on $\Delta^{\geq m+2}$:
\[ 0\to \cQ[e_{m}-e_{m+1}]\to  \cF^0/x_{m+1}\cF^0 \to \cC \to 0.\]
On the sheaf $\cF^0$ we can define a bilinear pairing by the same
formula as on $\cC$. This pairing is degenerate, but it induces the
pairing on $\cC$. Now as in Lemma~\ref{lem-exist} we construct a
homomorphism  $L: \cF^0\to\cF^0$ of degree $e_{m+1}$. We claim that this
homomorphism induces the injective homomorphism $\cC^0\to
\cC^1$. Indeed, we are reduced to the same Lemma~\ref{lem-last}. The
difference now is that we may have a strict inclusion $K\subset
K^\perp$, while the two spaces were equal in the proof of
Lemma~\ref{lem-exist}.
\end{remark}

Let us now put everything together and finish the proof of
Theorem~\ref{thm-main}. We start with the \Po sheaf
$\cL|_{\Delta^{\geq 1}}$ and apply the main construction to produce
new \Po sheaves $\cC$ and $\cQ$. Then inductively we apply the main
construction to $\cC$ and $\cQ$. These constructions give a Lefschetz
operation on each stalk $\cL_\sigma/x_d\cL_\sigma$,
$\dim\sigma=d$. Considering 
\[ \cL_{\hat{1}}/x_{n+1} \cL_{\hat{1}} \isom \Gamma(\cL,\Delta), \]
we get a Lefschetz operation on $\Gamma(\cL,\Delta)$ as stated in the
theorem.


\begin{thebibliography}{CC}

\bibitem{BBFK} G. Barthel, J.-P. Brasselet, K.-H. Fieseler and L. Kaup,
{\em Combinatorial Intersection cohomology for Fans}, T\^ohoku
Math. J. 54 (2002) 1-41.

\bibitem{BK} M.M. Bayer and A. Klapper,
{\em A new index for polytopes}, Discrete Comput. Geom. 6 (1991), 33-47.


\bibitem{BL} P. Bressler and V. A. Lunts, {\em Intersection cohomology
  on nonrational polytopes},  Compositio Math.  135  (2003),  no. 3,
  245-278.

\bibitem{B} M. Brion, {\em The Structure of the Polytope Algebra},
T\^ohoku Math. J. 49, 1997, 1-32.

\bibitem{F} W. Fulton, {\em Introduction to toric varieties},
  Princeton University Press, (1993).

\bibitem{K} K. Karu {\em The cd-index of fans and lattices}, preprint
  math.AG/0410513. 

\bibitem{S} R.P. Stanley, {\em Flag f-vectors and the cd-index},
  Math. Z. 216 (1994), 483-499.


\end{thebibliography}
\end{document}